\documentclass[12pt]{article}
\usepackage{amssymb,amsmath, amsthm, amscd,ifthen}
\usepackage[T2A]{fontenc}
\usepackage[cp1251]{inputenc}
\usepackage[english,russian]{babel}
\usepackage[dvips]{graphicx}
\usepackage{indentfirst}

\newtheorem{theorem}{Теорема}
\newtheorem{claim}{Утверждение}
\newtheorem{lemma}{Лемма}

\newtheorem{note}{Замечание}
\newtheorem{conjecture}{Гипотеза}
\newtheorem{example}{Пример}
\newtheorem{problem}{Проблема}

\theoremstyle{definition}

\newtheorem{defin}{Определение}

\begin{document}
\title{ Разложение расстановок чисел в кубе 
}
\author{И. Решетников}

\date{}

\maketitle

\begin{abstract}

Подмножество $M \subset \textbf{R}^3$ называется {\it базисным}, если для любой функции $f \colon M \to \textbf{R}$ существуют такие функции
$f_1; f_2; f_3 \colon \textbf{R} \to \textbf{R}$, что
$f(x_1, x_2, x_3) = f_1(x_1) + f_2(x_2) + f_3(x_3)$
для любой точки $(x_1, x_2, x_3)\in M$.
В работе доказан критерий базисности некоторых множеств в терминах   интересного свойства графов.
Приведены кон\-стру\-кции минимальных по включению небазисных множеств.


\end{abstract}

\section{Формулировки результатов}

\subsection{Постановка задачи и критерии базисности}

Понятие \textit{базисности} множества возникает в связи с тринадцатой проблемой Гильберта о суперпозиции непрерывных функций \cite{hilbert}, \cite{St}.

\begin{defin} \label{def_basic}
\footnote{
Другая формулировка базисности множества.
Андрей Николаевич и Владимир Игоревич играют в 3D-игру.
В кубе $n\times n\times n$, разбитом на $n^3$ единичных кубиков, отмечено
несколько кубиков.
А. Н. рас\-ста\-вля\-ет числа в отмеченных кубиках, как хочет.
В. И. смотрит на расставленные числа и выбирает $3n$ чисел $a_1,\dots,a_n,$ $b_1,\dots,b_n,$ $c_1,\dots,c_n$ 'весов слоёв'.
Если число в каждом отмеченном кубике $(i,j,k)$ (поставленное А. Н.) оказалось
равным сумме $a_i+b_j+c_k$ трёх 'весов слоёв',
содержащих этот кубик, то выиграл В. И.,
а иначе (т.е. если число хотя бы в одном отмеченном кубике оказалось не равным
сумме трех весов) выиграл А. Н.
Если выигрывает Владимир Игоревич, то множество центров отмеченных кубиков
базисное, в противном случае ---
небазисное.
}
Множество точек $M$ в $\textbf{R}^3$ называется
\emph{базисным}, если для любой функции $f \colon M \to \textbf{R}$ существуют такие функции $f_1, f_2, f_3 \colon \textbf{R} \to \textbf{R}$,
что $f( x_1, x_2, x_3 ) = f_1( x_1 ) + f_2( x_2 ) + f_3 ( x_3 )$ для любой точки $(x_1, x_2, x_3)\in M$. 

\end{defin}

    В других работах базисностью называется другое понятие. То, что мы называем ба\-зис\-нос\-тью, разумно называть разрывной базисностью \cite{hilbert}, но так как в данной работе используется только определение \ref{def_basic}, то мы опускаем слово 'разрывно'.  

Аналогично определяется базисность множеств точек в $\textbf{R}^2$ (см. подробнее в \cite{hilbert}, \cite{basic}).

Известен простой критерий
базисности множеств в $\textbf{R}^2$ \cite{hilbert}, \cite{basic}.

Легко понять, что базисность множества равносильна разрешимости некоторой системы ли\-ней\-ных уравнений.
Значит, существует алгоритм проверки базисности конечного мно\-же\-ства.

В этой работе мы приводим предварительные результаты, мотивированные следующей проблемой.

\begin{problem} \label{problem1} Найти 'простой' критерий
базисности конечного множества в \emph{$ \textbf{R}^3 $}.
\end{problem}


\begin{defin}

\emph{Слоем} будем называть плоскость, перпендикулярную оси координат.

\end{defin}

В частном случае, когда в каждом слое не более двух точек из множества $M$, базисность
множества
$M$ равносильна {\it базисности} некоторого графа.
Это позволяет найти простой критерий базисности множества $M$ для этого случая.

\begin{defin}[Базисность графа]
\footnote{Другая формулировка базисности графа.
Дан граф $G$. Два игрока играют в следующую игру.
Сначала первый игрок расставляет вещественные числа в вершинах графа $G$. Потом второй игрок расставляет вещественные числа на ребрах графа $G$. Действия игроков заканчиваются. Если число в каждой вершине равно сумме чисел на рёбрах, инцидентных ей, то выигрывает второй, иначе выигрывает первый.
Если выигрывает второй, то граф базисный.
}
Граф $G$ называется \emph{базисным}, если для любой расстановки вещественных чисел в вершинах графа $G$ существует расстановка ве\-щес\-твен\-ных чисел на рёбрах графа $G$ такая, что число в каждой вершине равно сумме чисел на рёбрах, инцидентных ей.
\end{defin}


\begin{theorem}\label{thKog}
Граф $G$ является базисным тогда и только тогда, когда
 в $G$ нет двудольной компоненты связности.
\end{theorem}

Доказательство теоремы \ref{thKog}, а также леммы \ref{lemma}, утверждений \ref{construction} и \ref{thBoyarov}, см. в параграфе 2.

Эта теорема интересна и сама по себе, а также в связи с определением одномерных когомологий графа \cite{prasolov}.

\begin{conjecture}
Аналог теоремы \ref{thKog} верен и для бесконечных графов, с за\-ме\-ча\-нием~\ref{note1}.
\end{conjecture}



\begin{claim}\label{thEq}
Пусть множество точек $M$ такое, что в каждом слое либо две точки из
$M$, либо ни одной. Рассмотрим граф $G(M)$ (возможно, бесконечный), вершинами которого являются точки множества $M$, причем вершины соединены ребром, если они лежат в одном слое.
Тогда множество $M$ является базисным тогда и только тогда, когда граф $G(M)$ базисный.
\end{claim}


Используя теорему \ref{thKog} и утверждение \ref{thEq}, можно построить контрпримеры к обращению следующего очевидного утверждения, доказываемого аналогично задаче 4a из \cite{basic}.  Самые простые из них можно построить и непосредственно, см. пример \ref{example1} в параграфе 2.

\begin{claim}[Достаточное условие базисности] \label{nonbasity}
Множество $M$ точек пространства
является базисным, если для любого непустого подмножества множества $M$
су\-щес\-тву\-ет слой, содержащий только одну точку из этого подмножества.
\end{claim}

Аналогичным способом можно сформулировать критерий базисности множества в $\textbf{R}^2$, что упрощает критерий базисности, приведённый в \cite{hilbert}.



\begin{lemma}[Необходимое условие базисности]
\label{lemma}
Если конечное множество \emph{$M \subset \textbf{R}^3$} можно раскрасить в два цвета так, чтобы в каждом слое было одинаковое количество точек каждого цвета, то множество $M$ не является базисным.
\end{lemma}


Обобщением леммы \ref{lemma} служит следующее очевидное утверждение, являющееся пере\-фор\-му\-ли\-ров\-кой утверждения о разрешимости системы линейных уравнений.

\begin{claim}[Критерий базисности множества] \label{criterium}
 Подмножество \emph{$M \subset \textbf{R}^3$} не является базисным тогда и только тогда, когда его точкам можно сопоставить целые 'веса', не все равные нулю так, чтобы для любого слоя $\alpha$ сумма 'весов' точек в пересечении $M \cap \alpha$ равнялась нулю.
\end{claim}

Аналог этого критерия справедлив и для плоскости.
Однако для плоскости 'веса' можно взять равными $0,+1,-1$.
Этот факт дает и красивый критерий базисности, и быстрый алгоритм проверки базисности.
Поэтому разумно пытаться найти критерий ба\-зи\-сно\-сти конечного множества в пространстве
в виде аналога утверждения \ref{criterium}, в котором были бы сильные ограничения на веса.
Пример \ref{ex2} показывает, что веса не всегда можно взять равными $0,+1,-1$.

\begin{problem} \label{pr2}
Можно ли ограничить модуль 'весов' в утверждении \ref{criterium} некоторой кон\-стан\-той?
\end{problem}


\subsection{Необходимые условия базисности}







Интересно приводить примеры
не базисных множеств и обобщать их до необходимых условий базисности.

\begin{defin}[Молния]
Пусть $\alpha$ --- слой. Обозначим координатные оси так, чтобы плоскость $\alpha$ была перпендикулярна оси $z$. Обозначим через $x(a)$ и $y(a)$ соответственно $x$- и $y$-координаты точки $a$. Конечная
последовательность точек $a_1, \ldots , a_n  \in \alpha$ называется молнией,
если для каждого $i$ выполнено $a_i \neq a_{i+1}$,
и при этом $x(a_i) = x(a_{i+1})$ для чётных $i$ и $y(a_i) = y(a_{i+1})$ для
нечётных $i$. (Не обязательно все точки молнии различны.)
Молния ${a_1, \ldots , a_{2l+1}}$ называется замкнутой, если $a_1 = a_{2l+1}$.
\end{defin}

В других работах то, что мы называем молнией, называется конечной молнией, а в определении молнии отсутствует требование конечности.

В двумерном случае минимальными по включению не базисными множествами явля\-ют\-ся замкнутые молнии без повторений вершин, и только они. Замкнутые молнии в слое явля\-ют\-ся минимальными по включению не базисными множествами и в пространстве, но не ис\-чер\-пы\-ва\-ют всё множество минимальных по включению не базисных множеств.

Приведём три конструкции минимальных по включению не базисных множеств в про\-стран\-стве.

\begin{claim} \label{construction} Раскрасим замкнутую молнию $M$ в слое $\alpha$ в чёрный и белый цвета так, чтобы соседние точки были покрашены в разные цвета. Разделим множество $M$ на несколько непересекающихся подмножеств так, чтобы в каждом подмножестве было одинаковое число точек каждого цвета. Каждое подмножество перенесём параллельно на вектор, перпендикулярный слою $\alpha$. Полученное множество точек
не базисное.
\end{claim}

Утверждение \ref{construction} следует из леммы \ref{lemma}.

Но даже среди подмножеств вершин куба есть отличный от множеств, получаемых в утверждении \ref{construction}, пример
минимального по включению не базисного множества.

\begin{example} \label{ex2}
Подмножество вершин куба $M := \{(0,0,0), (0,1,0), (1,0,0), (0,0,1), (1,1,1)\}$ не является базисным.

Более наглядно:


Первый слой $(z = 0)$
$$
  \begin{matrix}
    \boxtimes & \square \\
    \boxtimes & \boxtimes \\
  \end{matrix}
$$

второй слой $(z = 1)$
$$
  \begin{matrix}
    \square & \boxtimes \\
    \boxtimes & \square \\
  \end{matrix}
$$

\end{example}

Данный пример обобщается следующим образом.


\begin{claim}{ \label{thBoyarov} Пусть есть
не базисное множество $M$, в котором есть точки $A$~и~$B$, совпадающие по двум координатам. Рассмотрим множество $M' = ( M \setminus \{B\}) \cup \{A'\} \cup \{B'\}$, где точки $A'$ и $B'$ получаются параллельным переносом точек $A$ и $B$, соответственно, на вектор, перпендикулярный отрезку $AB$.
Множество $M'$ также
небазисное.}
\end{claim}


При построении не базисных множеств способами, указанными в утверждениях \ref{construction} и \ref{thBoyarov}, всегда найдётся пара точек, две координаты которых совпадают. Однако это не является не\-об\-хо\-ди\-мым условием базисности.
Существует не базисное множество, в котором нет двух точек, две координаты которых совпадают.


\begin{example} \label{cube}
Следующее множество не содержит двух точек, две координаты которых совпадают, но не является базисным.
$$\{(3, 0, 0), (2, 1, 0), (1, 2, 0), (0, 3, 0), (0, 0, 1), (1, 1, 1), (2, 2, 1), (3, 3, 1)\}.$$

Более наглядно:

Первый слой $(z = 0)$
$$
  \begin{matrix}
    \boxtimes & \square & \square &\square\\
    \square & \blacksquare& \square &\square \\
    \square& \square& \blacksquare& \square\\
    \square& \square&\square &\boxtimes \\
  \end{matrix}
$$

второй слой $(z = 1)$
$$
  \begin{matrix}
    \square  & \square  &\square  &\blacksquare\\
    \square  & \square & \boxtimes&\square  \\
    \square & \boxtimes& \square & \square \\
    \blacksquare& \square & \square &\square  \\
  \end{matrix}
$$

\end{example}

Не базисность этого примера следует из леммы \ref{lemma}. В наглядном изображении приведена раскраска, удовлетворяющая лемме \ref{lemma}.

\section{Доказательства}





\subsection{Доказательство теоремы \ref{thKog} }

\begin{proof}[Доказательство части 'только тогда' теоремы \ref{thKog}]
Достаточно доказать, что дву\-доль\-ный граф не является базисным. Если число в каждой вершине равно сумме чисел на рёбрах, инцидентных ей, то суммы чисел на вершинах долей равны. Поэтому для расстановки чисел, в которой в первой вершине стоит число 1,
а в остальных вершинах стоит число 0, не существует такой расстановки чисел на ребрах,
что число в каждой вершине равно сумме чисел на рёбрах, инцидентных ей.
\end{proof}

\begin{proof}[Доказательство части 'тогда' теоремы \ref{thKog}]
Часть 'тогда' вытекает из сформулированных ниже утверждений \ref{claim1} и \ref{claim2}.
\end{proof}

\begin{defin}
Кограницей $\delta v$ вершины $v$ будем называть рассстановку чисел на рёбрах графа, определенную следующим образом: на ребрах, инцидентных ей, стоит число 1, а на каждом из остальных ребер стоит число 0.
\end{defin}

Похожие, но другие, понятия кограниц возникают при определении одномерных ко\-го\-мо\-ло\-гий графа \cite{prasolov}.



\begin{claim}\label{claim1} Граф с $e$ рёбрами является базисным $\Leftrightarrow$ 
кограницы вершин линейно независимы в линейном пространстве \emph{$\textbf{R}^{e}$}.

\end{claim}

\begin{claim}\label{claim2} Кограницы вершин линейно независимы в линейном пространстве \emph{$\textbf{R}^{e}$} $\Leftarrow$ в графе $G$ нет двудольной компоненты связности.
\end{claim}

\begin{note} У утверждения \ref{claim2} верно и обращение.
\end{note}

\begin{proof}[Доказательство утверждения \ref{claim1}]
Пусть в вершинах стоят числа $b_1, \ldots, b_n$, а на рёбрах~--- числа $x_1, \ldots, x_e$. Составим матрицу $A$ системы линейных уравнений относительнно $x_1, \ldots, x_e$, которая обнуляется, если граф базисный. 
Положим $a_{ij} = 1$, если вершина~$v_i$ инцидентна ребру~$e_j$, и $a_{ij} = 0$ в противном случае. Таким образом строки матрицы $A$ есть кограницы вершин. Для любого столбца действительных чисел \textbf{b} cистема $A\textbf{x}=\textbf{b}$, где
$$
\textbf{b} =
  \begin{pmatrix}
    b_1 \\
    \vdots  \\
    b_n  \\
  \end{pmatrix}
,
\textbf{x} =
  \begin{pmatrix}
    x_1 \\
    \vdots  \\
    x_e  \\
  \end{pmatrix}
,
$$
разрешима $\Leftrightarrow$ строки матрицы $A$ линейно независимы.
\end{proof}

\begin{proof}[Доказательство утверждения \ref{claim2}]
Пусть есть линейная зависимость между кограницами вершин. Обозначим через $\lambda_i$ коэффициент при $\delta v_i$ в ней. Если вершины $u$ и $v$ соединены ребром, то $\lambda_v = - \lambda_u$. Поэтому кограницы вершин из компоненты связности вершины $v$ разбиваются на два класса эквивалентности --- с коэффициентом $\lambda_v$ и с коэффициентом $-\lambda_v$ соответственно. При этом вершины из каждого класса смежны только с вершинами из другого. Значит, эта компонента связности двудольна.
\end{proof}

\begin{note} \label{note1}
Аналог утверждения \ref{claim2} неверен для бесконечных графов. Контрпримером служит бесконечная цепочка.
\end{note}

\subsection{Другие доказательства}

\begin{example} \label{example1}
Базисное подмножество вершин куба $M := \{(0,1,0), (1,0,0), (0,0,1), (1,1,1)\}$.

Более наглядно:

Первый слой $(z = 0)$
$$
  \begin{matrix}
    \boxtimes & \square \\
    \square & \boxtimes \\
  \end{matrix}
$$

второй слой $(z = 1)$
$$
  \begin{matrix}
    \square & \boxtimes \\
    \boxtimes & \square \\
  \end{matrix}
.$$
\end{example}

\begin{proof}[Доказательство базисности множества из примера \ref{example1}]
Легко заметить, что граф $G$ для множества точек из примера \ref{example1} изоморфен графу $K_4$. В $K_4$ существует нечётный цикл. Из утверждения \ref{thEq} и теоремы \ref{thKog} следует, что множество $M$ в примере \ref{example1} является базисным.
\end{proof}

\begin{proof}[Прямое доказательство базисности множества из примера \ref{example1}]
Обозначим через $a, b, c, d$ значения функции $f$ в точках множества $M$. Соответственно, $a = f (0, 1, 0)$, $b = f  (1, 0, 0)$, $c = f (0, 0, 1)$, а $d =  f  (1, 1, 1)$. Тогда искомые функции $f_1, f_2, f_3$ можно определить следущим образом.
$$f_1(0)=a/2-d/2,\quad f_1(1)=b/2-c/2,\quad f_2(0)=b/2+c/2, \quad f_2(1)=a/2+d/2,$$
$$ f_3(0)=0, \quad f_3(1)=-a/2-b/2+c/2+d/2.$$
\end{proof}


\begin{proof}[Доказательство леммы \ref{lemma}]
Для любых функций $f_1, f_2, f_3 \colon M\rightarrow \textbf{R}$ будет выполнено, что сумма значений функции $f_1$ на проекциях точек белого цвета на ось $x$ равна сумме значений на проекциях точек черного цвета на ось $x$. То же для функций $f_2$ и $f_3$. А значит, функцию $f \colon M \rightarrow \textbf{R}$, которая в одной точке из $M$ принимает значение 1, а в остальных 0, нельзя представить в виде суммы функций $f_1, f_2$ и $f_3$.
\end{proof}

\begin{proof}[Доказательство утверждения \ref{construction}]
Для приведённой раскраски выполняется условие лем\-мы~\ref{lemma}, а значит, полученное множество не базисное. Из утверждения \ref{nonbasity} следует его ми\-ни\-маль\-ность по включению в случае, когда в каждом слое, перпендикулярном $\alpha$ и со\-дер\-жа\-щем точки множества, лежит ровно две точки из молнии.
\end{proof}

\begin{proof}[Доказательство утверждения \ref{thBoyarov}] 
Так как множество $M$ не базисное, то существует функ\-ция $f_b \colon M \rightarrow \textbf{R}$, не представимая в виде
$f_b(x_1, x_2, x_3) = f_1(x_1) + f_2(x_2) + f_3(x_3)$.
Тогда функ\-ция $f_b' \colon M' \rightarrow \textbf{R}$, которая на $(M \cap M') \setminus \{A'\} \setminus \{B'\}$ равна $f_b$, $f_b'(A'): = f_b(A)$, $f_b'(B'): = f_b(B)$, не представима в виде $f_b'(x_1, x_2, x_3) = f_1(x_1) + f_2(x_2) + f_3(x_3)$.

\end{proof}


Благодарю А. Б. Скопенкова за полезные обсуждения и замечания, Н. Волкова за написание программы, находящей примеры небазисных множеств среди подмножеств узлов куба $2 \times 2 \times 2$, И. Боярова за доказательство утверждения \ref{thBoyarov}.

\end{document}